\documentclass{article}
\usepackage{geometry}
\usepackage{graphicx}	
\usepackage[cp1251]{inputenc}
\usepackage[english]{babel}
\usepackage{mathtools}
\usepackage{amsfonts,amssymb,mathrsfs,amscd,amsmath,amsthm}
\usepackage{verbatim}

\def\ig#1#2#3#4{\begin{figure}[!ht]\begin{center}%
\includegraphics[height=#2\textheight]{#1.eps}\caption{#4}\label{#3}%
\end{center}\end{figure}}

\def\thtext#1{
  \catcode`@=11
  \gdef\@thmcountersep{. #1}
  \catcode`@=12
}

\def\threst{
  \catcode`@=11
  \gdef\@thmcountersep{.}
  \catcode`@=12
}

 \catcode`@=11
 \def\.{.\spacefactor\@m}
 \catcode`@=12

\theoremstyle{plain}
\newtheorem{thm}{Theorem}[section]
\newtheorem{prop}{Proposition}[section]
\newtheorem{cor}[prop]{Corollary}

\theoremstyle{definition}
\newtheorem{rk}[prop]{Remark}

\newcommand{\hR}{\widehat{R}}
\newcommand{\hx}{\hat{x}}
\newcommand{\hy}{\hat{y}}
\newcommand{\hX}{\widehat{X}}
\newcommand{\hY}{\widehat{Y}}

\newcommand{\cM}{\mathcal{M}}
\newcommand{\cP}{\mathcal{P}}
\newcommand{\cR}{\mathcal{R}}

\newcommand{\N}{\mathbb{N}}
\newcommand{\R}{\mathbb{R}}

\newcommand{\dl}{\delta}
\newcommand{\D}{\Delta}

\renewcommand{\r}{\rho}
\newcommand{\s}{\sigma}

\newcommand{\diam}{\operatorname{diam}}
\newcommand{\dis}{\operatorname{dis}}
\newcommand{\gen}{{\operatorname{gen}}}

\newcommand{\opt}{{\operatorname{opt}}}

\newcommand{\0}{\emptyset}
\renewcommand{\:}{\colon}

\renewcommand{\ss}{\subset}
\newcommand{\x}{\times}

\begin{document}
\title{Local Structure of Gromov--Hausdorff Space, and Isometric Embeddings of Finite Metric Spaces into this Space}
\author{Alexander O. Ivanov, Stavros Iliadis, Alexey A. Tuzhilin}
\date{}
\maketitle

\begin{abstract}
We investigate the geometry of the family $\cM$ of isometry classes of compact metric spaces, endowed with the Gromov--Hausdorff metric. We show that sufficiently small neighborhoods of generic finite spaces in the subspace of all finite metric spaces with the same number of points are isometric to some neighborhoods in the space $\R^N_{\infty}$, i.e., in the space $\R^N$ with the norm $\|(x_1,\ldots,x_N)\|=\max_i|x_i|$. As a corollary, we get that each finite metric space can be isometrically embedded into $\cM$ in such a way that its image belongs to a subspace consisting of all finite metric spaces with the same number $k$ of points. If the initial space has $n$ points, then one can take $k$ as the least possible integer with $n\le k(k-1)/2$.
\end{abstract}

%%%%%%%%%%%%%%%%%%%%%%%%%%%%%%
\section{Introduction}
\markright{\thesection.~Introduction}
%%%%%%%%%%%%%%%%%%%%%%%%%%%%%%
By  $\cM$ we denote the space of all compact metric spaces (considered up to an isometry) endowed with the Gromov--Hausdorff metric. It is well-known that $\cM$ is linear connected, complete, separable, but not proper. In a recent paper~\cite{IvaNikTuz}, A.~Ivanov, N.~Nikolaeva, and A.~Tuzhilin have shown that $\cM$ is geodesic. There are many other open questions concerning geometrical properties of $\cM$. S.~Iliadis formulated the following problem: \emph{is it true that $\cM$ is universal for the family of compact metric spaces\/}? The latter may be interpreted in weak and strong senses. The weak sense means that any compact metric space can be isometrically embedded into $\cM$. The strong sense means in addition that each isometric embedding of a subspace of a compact metric space can be extended to an isometric embedding of the whole space. It is easy to see that $\cM$ is not universal in the strong sense (see below). Concerning the weak universality, we show the following: \emph{each finite metric space can be isometrically embedded into $\cM$}. Moreover, we construct such an embedding with its image belonging to the subspace of all finite metric spaces with $k$ points: if the initial space has $n$ points, then one can  chose $k$ as the least possible integer such that $n\le k(k-1)/2$.

The construction of such embedding is based on our results concerning the local geometry of the family of finite metric spaces with fixed number of points considered in sufficiently small neighborhoods of generic spaces. More precisely, we show that such neighborhoods are isometric to some neighborhoods of the corresponding points in the space $\R^k_{\infty}$, i.e., in the space $\R^k$ with the norm $\bigl\|(x_1,\ldots,x_k)\bigr\|=\max_i|x_i|$.

%%%%%%%%%%%%%%%%%%%%%%%%%%%%%%
\section{Preliminaries}
\markright{\thesection.~Preliminaries}
%%%%%%%%%%%%%%%%%%%%%%%%%%%%%%
Let $X$ be an arbitrary metric space.  By $|xy|$ we denote the distance between points $x$ and $y$ in $X$. For every point $x\in X$ and a real number $r>0$ by $U_r(x)$ we denote the open ball of radius $r$ centered at $x$; for every nonempty $A\ss X$ and real number $r>0$ we put $U_r(A)=\cup_{a\in A}U_r(a)$.

For nonempty $A,\,B\ss X$, let us put
$$
d_H(A,B)=\inf\bigl\{r>0:A\ss U_r(B)\&B\ss U_r(A)\bigr\}.
$$
This value is called the \emph{Hausdorff distance between $A$ and $B$}. It is well-known~\cite{BurBurIva} that the restriction of the Hausdorff distance to the family of all closed bounded subsets of $X$ is a metric.

Let $X$ and $Y$ be metric spaces. A triple $(X',Y',Z)$ that consists of a metric space $Z$ and its subsets $X'$ and $Y'$ isometric to $X$ and $Y$, respectively, is called a \emph{realization of the pair $(X,Y)$}. The \emph{Gromov--Hausdorff distance $d_{GH}(X,Y)$ between $X$ and $Y$} is the greatest lower bound of the real numbers $r$ such that there exists a realization $(X',Y',Z)$ of the pair $(X,Y)$ with $d_H(X',Y')\le r$. It is well-known~\cite{BurBurIva} that the $d_{GH}$ restricted to the family $\cM$ of isometry classes of compact metric spaces is a metric.

Recall that a \emph{relation\/} between sets $X$ and $Y$ is a subset of the Cartesian product $X\x Y$. By $\cP(X,Y)$ we denote the set of all nonempty relations between $X$ and $Y$. If $\pi_X\:(x,y)\mapsto x$ and $\pi_Y\:(x,y)\mapsto y$ are the canonical projections, then their restrictions to each $\s\in\cP(X,Y)$ are denoted in the same manner.

We consider each relation $\s\in\cP(X,Y)$ as a multivalued mapping, whose domain may be less than the whole $X$. By analogy with with mappings, for every $x\in X$ its image $\s(x)=\{y\in Y\mid(x,y)\in\s\}$ is defined, and for every $y\in Y$  its preimage $\s^{-1}(y)=\{x\in X\mid(x,y)\in\s\}$ is defined also; for every  $A\ss X$ its image $\s(A)$ is the union of the images of all the elements from $A$, and, similarly, for every $B\ss Y$ its preimage is the union of the preimages of all the  elements from $B$.

A relation $R$ between $X$ and $Y$ is called a \emph{correspondence}, if the restrictions of the canonical projections $\pi_X$ and $\pi_Y$ onto $R$ are surjections.  By $\cR(X,Y)$ we denote the set of all correspondences between $X$ and $Y$.

Let $X$ and $Y$ be metric spaces, then for every relation $\s\in\cP(X,Y)$ its \emph{distortion $\dis\s$} is defined as
$$
\dis\s=\sup\Bigl\{\bigl||xx'|-|yy'|\bigr|: (x,y)\in\s,\ (x',y')\in\s\Bigr\}.
$$

The following result is well-known.

\begin{prop}[\cite{BurBurIva}]\label{th:GH-metri-and-relations}
For any metric spaces $X$ and $Y$ we have
$$
d_{GH}(X,Y)=\frac12\inf\bigl\{\dis R\mid R\in\cR(X,Y)\bigr\}.
$$
\end{prop}

If $X$ and $Y$ are finite metric spaces, then the set $\cR(X,Y)$ is finite, hence, there exists an $R\in\cR(X,Y)$ such that $d_{GH}(X,Y)=\frac12\dis R$. Every such correspondence $R$ is called \emph{optimal}. Notice that optimal correspondences do also exist for any compact metric spaces $X$ and $Y$, see~\cite{IvaIliadisTuz}.  By $\cR_\opt(X,Y)$ we denote the set of all optimal correspondences between $X$ and $Y$.

For arbitrary nonempty sets $X$ and $Y$ a correspondence $R\in\cR(X,Y)$ is called \emph{irreducible}, if it is a minimal element of the set $\cR(X,Y)$ w.r.t. the order given by the inclusion relation.  By $\cR^0(X,Y)$  we denote the set of all nonempty irreducible correspondences between $X$ and $Y$.

The next result describes the structure of irreducible correspondences.

\begin{prop}\label{prop:IrreducibleRelation}
Every irreducible correspondence $R\in\cR(X,Y)$ generates partitions $X=X'_1\sqcup X''_1\sqcup X_2$ and $Y=Y'_1\sqcup Y''_1\sqcup Y_2$, together with the partitions $\hX'_1$ and $\hY'_1$ of the sets $X'_1$ and $Y'_1$, respectively. Also, $R$ induces a bijection $\hR$ between the sets $\hX=\hX'_1\sqcup X''_1\sqcup X_2$ and $\hY=Y_2\sqcup Y''_1\sqcup\hY'_1$ such that $(x,y)\in R$ iff either $x\in\hx\in\hX'_1,\,y=\hR(\hx)\in Y_2$, or $x\in X''_1,\,y=\hR(y)\in Y''_1$, or $x\in X_2,\,y\in\hy=\hR(x)\in\hY'_1$.
\end{prop}

Figure~\ref{fig:IrreducibleRelation} illustrated the latter theorem.

\ig{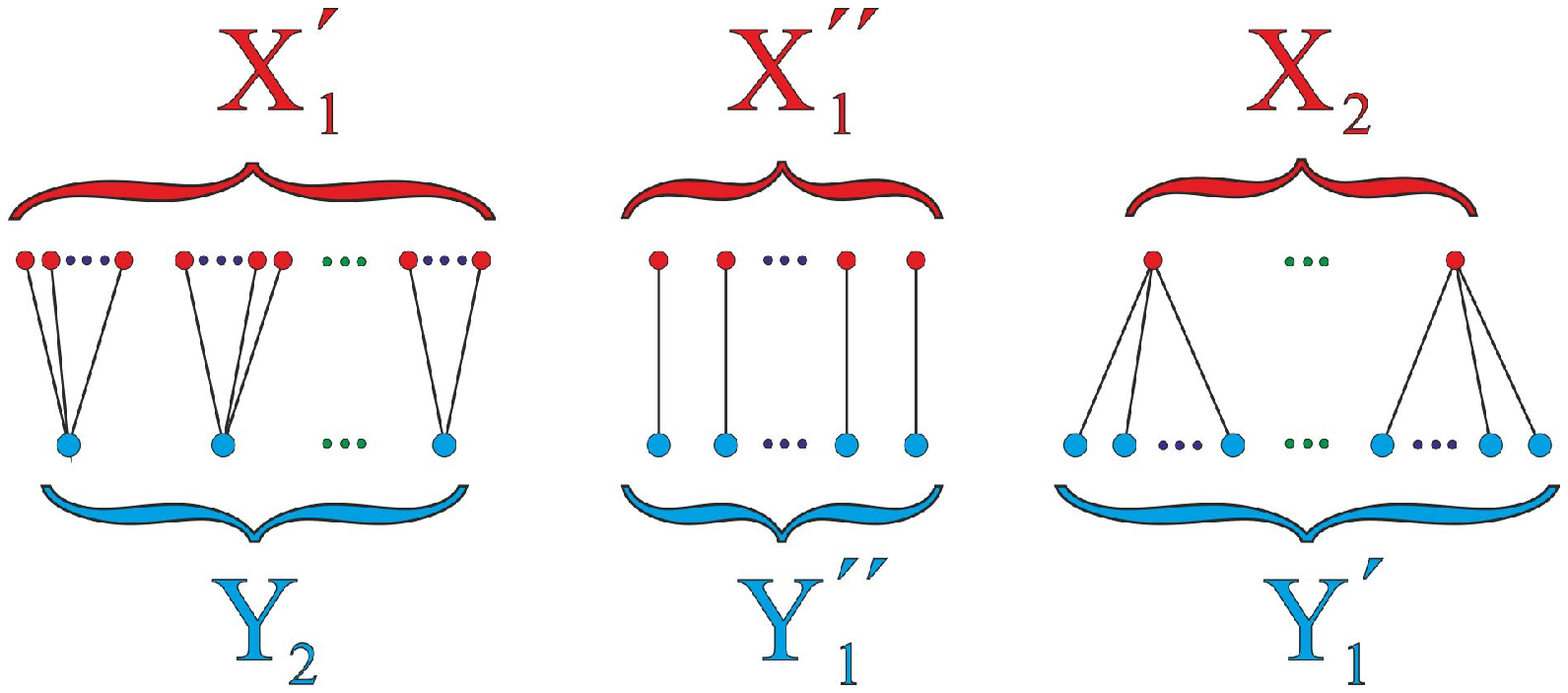}{0.25}{fig:IrreducibleRelation}{The structure of an irreducible correspondence.}

The next results (see also~\cite{IvaTuzIrreducible}) demonstrates  the importance of the irreducible correspondences for calculating the Gromov--Hausdorff distances.

\begin{prop}
For every $R\in\cR(X,Y)$ there exists an $R^0\in\cR^0(X,Y)$ such that $R^0\ss R$.
\end{prop}

Now, let $X$ and $Y$ be metric spaces. We put $\cR^0_\opt(X,Y)=\cR^0(X,Y)\cap\cR_\opt(X,Y)$.

\begin{cor}
For any $X,Y\in\cM$ we have $\cR^0_\opt(X,Y)\ne\0$.
\end{cor}

If $X$ is a metric space, then by $\diam X$ we denote the diameter of $X$, i.e., the value 
$$
\diam X=\sup\bigl\{|xx'|:x,x'\in X\bigr\}.
$$ 
The next result is well-known~\cite{BurBurIva}.

\begin{prop}\label{prop:dist_and_diams}\hfill\hbox{}
\begin{enumerate}
\item\label{prop:dist_and_diams:1} If $\D_1$ stands for the one-point metric space, and $X$ is an arbitrary metric space, then
$$
d_{GH}(\D_1,X)=\frac12\diam X.
$$
\item\label{prop:dist_and_diams:2} For any metric spaces $X$ and $Y$ it holds
$$
d_{GH}(X,Y)\le\frac12\max\{\diam X,\,\diam Y\}.
$$
\end{enumerate}
\end{prop}

%%%%%%%%%%%%%%%%%%%%%%%%%%%%%%
\section{Strong Non-Universality of the Gromov--Hausdorff Space}
\markright{\thesection.~Strong nonuniversality of the Gromov--Hausdorff space}
%%%%%%%%%%%%%%%%%%%%%%%%%%%%%%

Show that the two-point metric space $\{A,B\}$ with the distance $|AB|=1/2$ may be isometrically embedded into $\cM$ in such a way that this embedding can not be extended to an isometric embedding of three-point metric space $\{A,B,C\}$ with distances $|AC|=1/2$ and $|BC|=2/3$. This example demonstrates that the space $\cM$ is not universal in the strong sense.

Indeed, we map the $A$ to the one-point space $\D_1$, and the $B$ to the two-point space $\D_2$ with distance $1$. By \ref{prop:dist_and_diams}~(\ref{prop:dist_and_diams:1}), we have $d_{GH}(\D_1,\D_2)=\frac12\diam\D_2=1/2$, thus, we have embedded isometrically the $\{A,B\}$ into $\cM$.

Now, suppose that the $C$ is mapped to a point $X\in\cM$ such that $|AC|=d_{GH}(\D_1,X)$. By \ref{prop:dist_and_diams}~(\ref{prop:dist_and_diams:1}), $|AC|=\frac12\diam X=1/2$, hence, $\diam X=1$. However, by \ref{prop:dist_and_diams}~(\ref{prop:dist_and_diams:2}), $d_{GH}(\D_2,X)\le\frac12\max\{\diam\D_2,\,\diam X\}=1/2$. Thus, the distance  $|X\D_2|$ can not be equal to $2/3$ for any choice of the $X$.

%%%%%%%%%%%%%%%%%%%%%%%%%%%%%%
\section{Geometry of Gromov--Hausdorff Space in Neighborhoods of Generic Spaces}
\markright{\thesection.~Local structure of Gromov--Hausdorff space}
%%%%%%%%%%%%%%%%%%%%%%%%%%%%%%

By $\cM_n$ we denote  the subset of $\cM$ consisting of all metric spaces each of which has at most $n$ points; let $\cM_{[n]}$ stand for the subset of $\cM_n$ consisting of all $n$-point metric spaces.

For every $n\in\N$ we define a mapping $\nu\:\cM_{[n]}\to\R^{n(n-1)/2}_{\infty}$ in the following way. For $X=\{x_1,\ldots,x_n\}\in\cM_{[n]}$ we consider all the distances $|x_ix_j|$ between distinct points, arrange these values in ascending order, and  let $\nu(X)$ be the resulting vector from $\R^{n(n-1)/2}$ divided by $2$. Let $\nu^i(X)$ be the $i$th coordinate of the vector $\nu(X)$.

We say that $X,Y\in\cM_{[n]}$ are \emph{structural isomorphic\/} if there exists a bijection $f\:X\to Y$ preserving the order on the set of distances: $|xy|\le |zw|$ iff $|f(x)f(y)|\le|f(z)f(w)|$. Every such $f$ we call a \emph{structural isomorphism}. Notice that for a structural isomorphism $f\:X\to Y$ and each $1\le i\le n(n-1)/2$ the condition $\nu^i(X)=|xx'|/2$ implies $\nu^i(Y)=\bigl|f(x)f(x')\bigr|/2$. In other words, for any $i$, the $i$th coordinates of vectors $\nu(X)$ and $\nu(Y)$ are equal to half-distances between the pairs of points from $X$ and $Y$ corresponding to each other under the structural isomorphism $f$. Thus, if we consider $f$ as a correspondence between $X$ and $Y$, i.e., $f\in\cR(X,Y)$, then we get $\big|\nu(X)\nu(Y)\big|=\frac12\dis f$.

By $\cM'_{[n]}$ we denote the subset of $\cM_{[n]}$ consisting of all the spaces such that all non-zero distances in them are distinct. Notice that if $X,Y\in\cM'_{[n]}$ are structural isomorphic, then structural isomorphism is uniquely defined. Moreover, the structural isomorphism generates an equivalence relation on $\cM'_{[n]}$, therefore, $\cM'_{[n]}$ is partitioned into the corresponding equivalence classes that are referred as \emph{classes of structural isomorphism}. Notice the following obvious property: for any structural isomorphic $X,Y\in\cM'_{[n]}$, a space $Z\in\cM_{[n]}$ is structurally isomorphic to $X$ iff it is structurally isomorphic to $Y$. The latter proves the correctness of the concept of the \emph{closure of a structural isomorphism class\/} as of the family of all $Z\in\cM_{[n]}$ which are structurally isomorphic to some $X\in\cM'_{[n]}$ (this definition does not depend on the choice of the representative $X$ in the structural isomorphism class). Thus, the space $\cM_{[n]}$ is covered  by the closures of the structural isomorphism classes, and each two spaces from the same such closure are structurally isomorphic. It is not difficult to show that these closures can be defined as maximal subfamilies of $\cM_{[n]}$ such that any two their elements are structurally isomorphic.

The next result is not used in the present paper but seems to be self-important.

\begin{prop}
Let $C\ss\cM_{[n]}$ be the closure of a structural isomorphism class. Then the mapping $\nu\:C\to\R^{n(n-1)/2}_{\infty}$ is incompressible.
\end{prop}

\begin{proof}
Let us choose arbitrary $X,Y\in C$, and let $R\:X\to Y$ be a structural isomorphism. Then $R\in\cR(X,Y)$, therefore,
$$
d_{GH}(X,Y)\le\frac12\dis R=|\nu(X)\nu(Y)|.
$$
\end{proof}

By a \emph{generic space\/} we mean every finite metric space such that all its non-zero distances are pairwise distinct, and all its triangle inequality hold strictly.   By $\cM^\gen$ we denote the family of all the generic spaces. Clearly that $\cM^\gen$ is everywhere dense in $\cM$. We put $\cM_{[n]}^\gen=\cM_{[n]}\cap\cM^\gen$. By definition, $\cM_{[n]}^\gen\ss\cM'_{[n]}$, thus, $\cM_{[n]}^\gen$ is also partitioned into the structural isomorphism classes.

Let $X=\{x_1,\ldots,x_n\}\in\cM$, $n\ge3$. We define $\dl(X)$ to be equal to the least of the following two numbers:
\begin{flalign*}
\indent&\min\bigl\{|x_ix_j|+|x_jx_k|-|x_ix_k|:\#\{i,j,k\}=3\bigr\},&\\
\indent&\min\Bigl\{\bigl||x_ix_j|-|x_px_q|\bigr|:\#\{i,j,p,q\}\ge3\Bigr\}.&
\end{flalign*}
For $n=2$ we put $\dl(X)=|x_1x_2|$.

\begin{rk}
If we put $p=q$ in the definition of the second minimum, then we get $|x_ix_j|$. Thus, $\dl(X)\le|x_ix_j|$ for all $i\ne j$. For $n=2$ this property holds as well.
\end{rk}

It is easy to see that $X=\{x_1,\ldots,x_m\}$ belongs to $\cM^\gen$ iff $\dl(X)>0$. Besides, for any $Y=\{y_1,\ldots,y_n\}\in\cM$ such that $\bigl||x_ix_j|-|y_iy_j|\bigr|<\dl(X)/3$ for all $1\le i,\,j\le n$ we have $Y\in\cM^\gen$, the $X$ and $Y$ are structurally isomorphic, and the mapping $x_i\mapsto y_i$ is a structural isomorphism.

\begin{thm}\label{thm:isometry}
Let $X=\{x_1,\ldots,x_n\}\in\cM^\gen$. We put $\dl=\dl(X)/6$, $U=U_\dl(X)\ss\cM_n$, $N=n(n-1)/2$, $V=U_\dl\bigl(\nu(X)\bigr)\ss\R^N_\infty$. Then the mapping $\nu|_U\:U\to V$ is an isometry.
\end{thm}

\begin{proof}
Choose an arbitrary $Y\in U$. Let us show that $\#Y=n$. Indeed, suppose otherwise that $\#Y<n$, then for any $R\in\cR(X,Y)$ there exists an $y\in Y$ such that for some distinct $x_1,x_2\in X$ it holds $(x_1,y),\,(x_2,y)\in R$, hence, $\dis R\ge|x_1x_2|\ge\dl(X)=6\dl$. Therefore, $d_{GH}(X,Y)\ge3\dl$, a contradiction.

So, $Y=\{y_1,\ldots,y_n\}$. Choose an arbitrary $R\in\cR^0_\opt(X,Y)$. If $R$ is not a bijection, then we can apply the arguments we used just above and come to a contradiction again. Thus, $R$ has to be a bijection. Renumbering if necessary the elements from $Y$, without loss of generality we can assume that $y_i=R(x_i)$ for all $1\le i\le n$. Now, show that $R$ is a structural isomorphism. Notice that for any $1\le i,\,j\le n$ the relations
$$
\bigl||y_iy_j|-|x_ix_j|\bigr|\le\dis R=2d_{GH}(X,Y)<2\dl=\dl(X)/3
$$
are valid, and hence, if $0<|x_ix_j|<|x_px_q|$, then we have $|y_py_q|-|y_iy_j|>\dl(X)/3>0$, because $|x_px_q|-|x_ix_j|\ge\dl(X)$.

Therefore, $\big|\nu(X)\nu(Y)\big|=\frac12\dis R=d_{GH}(X,Y)$, thus, the mapping $\nu$ is an isometry.

It remains  to show that $\nu$ is surjective. Choose an arbitrary $Z=(\r_{12},\ldots,\r_{(n-1)n})\in V$, and for every $1\le i<j\le n$ put $\r_{ji}=\r_{ij}$. Since $Z\in U_\dl\bigl(\nu(X)\bigr)$, we have $\bigl|\r_{ij}-|x_ix_j|/2\bigr|<\dl$ for all distinct $1\le i,\,j\le n$. Choose arbitrary pairwise distinct $1\le i,\,j,\,k\le n$, then
$$
\r_{ij}+\r_{jk}-\r_{ik}>\frac12\bigl(|x_ix_j|+|x_jx_k|-|x_ix_k|\bigr)-3\dl\ge\dl(X)/2-\dl(X)/2=0.
$$
Besides  $\r_{ij}>|x_ix_j|/2-\dl\ge\dl(X)/2-\dl(X)/6>0$. Thus, the values $2\r_{ij}$ generates a metric $|y_iy_j|=2\r_{ij}$ on the set $\{y_1,\ldots,y_n\}$. Let $Y=\{y_1,\ldots,y_n\}$ be the metric space obtained in this way, and let $R\in\cR(X,Y)$ be the bijection $x_i\mapsto y_i$. Then
$$
d_{GH}(X,Y)\le\frac12\dis R=\frac12\max\Bigl\{\bigl||y_iy_j|-|x_ix_j|\bigr|\Bigr\}= \max\Bigl\{\bigl|\r_{ij}-|x_ix_j|/2\bigr|\Bigr\}<\dl,
$$
thus, $Y\in U_\dl(X)=U$. It remains to notice that $\nu(Y)=Z$.
\end{proof}

%%%%%%%%%%%%%%%%%%%%%%%%%%%%%%
\section{Isometric Embedding of a Finite Metric Space into $\cM$}
\markright{\thesection.~Isometric embedding of a finite metric space into $\cM$}
%%%%%%%%%%%%%%%%%%%%%%%%%%%%%%

In this section we prove the weak universality of the Gromov--Hausdorff space for finite metric spaces.

\begin{thm}
Let $X$ be an arbitrary finite metric space consisting of $n$ points, and $k$ be the least integer such that $n\le k(k-1)/2$. Then $X$ can be isometrically embedded into $\cM$ in such a way that its image belongs to $\cM_{[k]}$.
\end{thm}

\begin{proof}
Put $X=\{x_1,\ldots,x_n\}$. First assume that $n=k(k-1)/2$ for some $k\in\N$.

By $f\:X\to\R^n_\infty$ we denote the isometric Kuratowski embedding: $f\:x_i\mapsto\bigl(|x_1x_i|,\ldots,|x_nx_i|\bigr)$, see~\cite{ITMinFil}. Recall that the translations  are isometries of $\R^n_\infty$.

Let $d$ stand for the diameter of $X$. Consider a generic metric space $S\in\cM_{[k]}$ such that $\dl(S)/6>d$. Put $\dl=\dl(S)/6$, then $Z=\nu(S)+f(X)\ss U_\dl\bigl(\nu(S)\bigr)\ss\R^n_\infty$. As we mentioned above, the space $Z$ with the metric induced from $\R^n_\infty$ is isometric to $X$.

Put $U=U_\dl(S)\ss\cM_{[k]}$ and $V=U_\dl(\nu(S))\ss\R^n_\infty$, then, by Theorem~\ref{thm:isometry}, $\nu|_U\:U\to V$ is an isometry. Therefore, $\nu^{-1}(Z)\ss\cM_{[k]}\ss\cM$ is isometric to $X$.

Now, assume that $n\ne k(k-1)/2$ for any $k\in\N$. Consider the least $k$ such that $n<k(k-1)/2$. Extend $X$ upto a metric space $Y$ consisting of $k(k-1)/2$ points. We do it in such a way that the distances between points from $X$ are preserved, and all the remaining distances are equal to the diameter $d$ of the space $X$ (it is easy to verify that $Y$ is a metric space). By the previous arguments, $Y$ is isometric to a subspace of $\cM_{[k]}$, thus, $X$ is isometric to the corresponding part of this subspace.
\end{proof}

\end{document}